\documentclass[12pt,reqno]{amsart}
\usepackage{fullpage}
\usepackage{amsfonts}
\usepackage{mathpazo}
\usepackage{graphicx}

\newcommand{\Q}{{\mathbb Q}}
\newcommand{\R}{{\mathbb R}}

\renewcommand{\L}{{\mathcal L}}
\renewcommand{\S}{{\mathcal S}}
\newcommand{\css}{circle-center set}
\newcommand{\vvec}[1]{\overrightarrow{#1}}
\newcommand{\ep}{\varepsilon}

\newtheorem{problem}{Problem}

\newtheorem{exercise}{Exercise}
\newtheorem*{definition}{Definition}
\newtheorem{theorem}{Theorem}
{
\theoremstyle{remark}
\newtheorem{example}{Example}
\newtheorem*{remark}{Remark}
}

\begin{document}

\title{Sets that contain their circle centers}
\author{Greg Martin \\ \today}
\address{Department of Mathematics \\ University of British Columbia \\ Room
121, 1984 Mathematics Road \\ Canada V6T 1Z2}
\email{gerg@math.ubc.ca}
\maketitle

\section{Introduction}

In 2001, Stan Wagon posed an interesting problem on the Macalester College Problem of the Week website \cite{PoW}:

\begin{problem}
Let $\S$ be a finite set of points in the plane in general position, that is, no three points of $\S$ lie on a line. Show that if there are at least three points in $\S$, then $\S$ cannot ``contain its circle centers''.
\end{problem}

The terminology needs a brief explanation. A set $\S\subset\R^2$ is said to {\it contain its circle centers\/} if, for any three non-collinear points from $\S$, the center of the circle through those points is always in $\S$; in other words, if $S$ contains the vertices of any triangle, then it also contains the triangle's circumcenter.

Several solutions were quickly submitted, and two (outlined in Exercises \ref{minex} and \ref{iterex} below) were posted on the website, along with some discussion. An especially interesting feature of this pair of solutions was that one did not use the assumption that no three points of $\S$ lie on a line, while the other did not use the assumption that $\S$ is finite!

The discussion naturally turned at that point to how much could be said about sets that contain their circle centers. The entire plane $\S=\R^2$ is a perfectly reasonable example of such a set; less trivially, the set $\S=\Q^2$, consisting of all points in the plane with rational numbers as coordinates, contains its circle centers (see Exercise \ref{q2ex} below). Are there lots of sets with this property, or only a few?

It's clearly cheating to put {\em all} the points of $\S$ on a single line, or else there are no circle centers to form. This motivates the following definition:

\begin{definition}
A {\em\css} is a subset of $\R^2$ that is not a subset of any line and that contains its circle centers.
\end{definition}

John Guilford conjectured that any \css\ must be unbounded. Working on this conjecture eventually led to an even stronger discovery: {\em there is essentially only one \css}. The reader might well be skeptical at this point, especially given the fact that we have already mentioned {\em two} \css s, namely $\R^2$ and $\Q^2$!

To pin down what we mean by ``essentially'', we need to review some basic topological terms about sets in the plane. An {\em open disk} is simply the interior of any circle. A set $\S$ is {\em dense} if every open disk contains at least one of the set's points, or equivalently if for any point $P$ in $\R^2$, there is a sequence of points $\{P_1,P_2,P_3,\dots\}$ from $\S$ that converges to~$P$. A {\em closed set} is a set $\S$ with the property that, for any sequence of points $\{P_1,P_2,P_3,\dots\}$ from $\S$ whose limit exists, the limit itself must be a point in~$\S$. Finally, the {\em closure} $\overline\S$ of $\S$ is the set of all points we can obtain by taking limits of sequences of points from $\S$, or equivalently the smallest closed set containing $\S$. (See \cite[Chapter 1]{DR} for information about these topological concepts in more abstract spaces.)

Using this terminology, the first version of our theorem is:

\begin{theorem}
Every \css\ $\S$ must be dense in $\R^2$.
\label{dense.thm}
\end{theorem}

The \css\ $\Q^2$ is indeed dense in $\R^2$: every open disk contains a point both of whose coordinates are rational.

We can rephrase this theorem slightly by pointing out that if a set $\S$ contains its circle centers, then its closure $\overline\S$ also contains its circle centers (see Exercise \ref{closureex} below). Since the closure of any dense subset of $\R^2$ is $\R^2$ itself, the theorem above is actually equivalent to the following version:

\begin{theorem}
$\R^2$ is the only closed \css.
\label{only1.thm}
\end{theorem}

This is the sense in which there is ``only one'' set containing its circle centers. The example $\Q^2$ does not contradict this statement, since it is not closed, and indeed the closure of $\Q^2$ is all of $\R^2$.

In the following exercises and the rest of this paper, we use the notation $C(P,Q,R)$ to denote the circumcenter of $\triangle PQR$ and $r(P,Q,R)$ to denote its circumradius, that is, $C(P,Q,R)$ and $r(P,Q,R)$ are the center and radius, respectively, of the circle going through the points $P$, $Q$, and~$R$.

\begin{exercise}\label{minex}
Given three non-collinear points $D,E,F$, let $G = C(D,E,F)$ and $H=C(D,E,G)$. Show that one of $r(D,E,G)$, $r(D,F,G)$, and $r(E,F,G)$ is less than $r(D,E,F)$, unless $\triangle DEF$ is equilateral, in which case $r(D,G,H) < r(D,E,F)$. Conclude that no \css\ can be finite by considering, for a potential counterexample $\S$, the minimal radius of a circle through three points of $\S$.
\end{exercise}

\begin{exercise}\label{iterex}
Given three non-collinear points $D,E,F$, let $G = C(D,E,F)$, $H=C(D,E,G)$, and $I=C(D,E,H)$. Show that either one of $G$, $H$, or $I$ lies on the line $DE$, or else the three points $G$, $H$, and $I$ are collinear. Furthermore, show that $G$, $H$, and $I$ are not distinct, then $\triangle DGH$ is equilateral and $C(D,G,H)$ lies on the line $DE$. Conclude that no \css\ can be in general position.
\end{exercise}

\begin{exercise}\label{q2ex}
Suppose that both endpoints of a line segment are in $\Q^2$. Show that the perpendicular bisector of the line segment has an equation that can be written in the form $ax+by+c=0$ with $a$, $b$, and $c$ rational. Conclude that if $D$, $E$, and $F$ are non-collinear points in $\Q^2$, then the center of the circle through $D$, $E$, and $F$ is also in $\Q^2$.
\end{exercise}

\begin{exercise}\label{closureex}
Show that the location of the circumcenter of a triangle depends continuously on the positions of the triangle's vertices. Conclude that the circumcenter $C(\lim P_n,\lim Q_n,\lim R_n)$ of three limit points is the same as the limiting circumcenter $\lim C(P_n,Q_n,R_n)$, and so the closure of a \css\ is again a \css.
\end{exercise}

\section{First dimension first}
\label{cotsec}

Before we attack the problem stated in the introduction, let's first consider a ``one-dimensional'' version of the problem. Note that every \css\ contains isosceles triangles: if $M$, $N$, and $P$ are points in a \css, then $Q=C(M,N,P)$ is also in the \css, and each of $\triangle QMN$, $\triangle QNP$, and $\triangle QPM$ is isosceles. We can scale and rotate triangles freely without changing the relative positions of their circumcenters, so we don't lose generality by looking at the specific case of an isosceles triangle whose base runs vertically from $M=(0,1)$ to $N=(0,-1)$.

In this section we consider only circumcenters of isosceles triangles with $M$ and $N$ as vertices. Note that the third vertices of these triangles, and their circumcenters, all lie on the $x$-axis, and so this really is a ``one-dimensional'' version of the original circle-center problem. So let's look closely at how the circumcenter of an isosceles triangle relates to the dimensions of that triangle.
 
\newtheorem*{CT}{Cotangent Relation}
\begin{CT}
Let $\alpha$ be a real number that is not a multiple of that is not a multiple of $\pi/2$.  The center of the circle passing through the three points $(0,1)$, $(0,-1)$, and $(\cot\alpha,0)$ is $(\cot2\alpha,0)$.
\end{CT}

\begin{proof}[Proof when $0<\alpha<\pi/4$]
Define the points $M=(0,1)$, $N=(0,-1)$, $O=(0,0)$, and $P=(1,1)$. Let $Q=(\cot\alpha,0)$, so that the measure of $\angle MQO$ equals $\alpha$. Since $\overline{OQ}$ and $\overline{MP}$ are parallel, we see that the measure of $\angle QMP$ equals $\alpha$ as well (see Figure~\ref{CRfig}). If we now let $R=(\cot2\alpha,0)$, then the same argument shows that both $\angle MRO$ and $\angle RMP$ have measure $2\alpha$. However, this implies (as Figure~\ref{CRfig} shows) that $\angle MQR$ and $\angle QMR$ both have measure $\alpha$, and so $\triangle MQR$ is isosceles with $QR=MR$.

%:CRfig
\begin{figure}[hbtp]

\begin{picture}(0,0)
\put(-4, 104){$M$}
\put(-3, 2){$O$}
\put(120, 104){$P$}
\put(373, 6){$Q$}
\put(177, 6){$R$}
\put(332, 6){$\scriptstyle\alpha$}
\put(44, 88.5){$\scriptstyle\alpha$}
\put(49, 98.5){$\scriptstyle\alpha$}
\put(150, 6){$\scriptstyle2\alpha$}
\end{picture}
\includegraphics[height=1.5in]{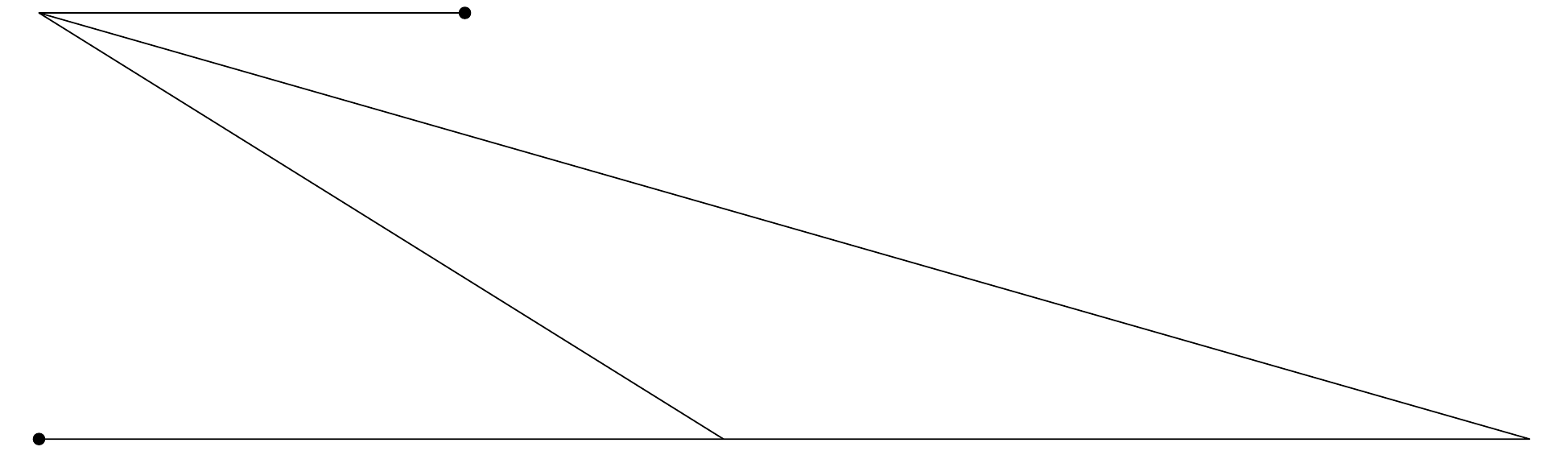}

\caption{The Cotangent Relation}
\label{CRfig}
\end{figure}
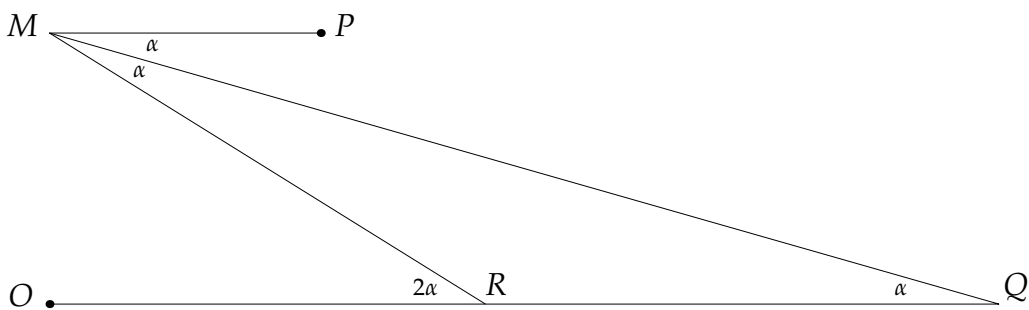

A symmetric argument shows that $QR=NR$ as well. Since $R$ is equidistant from the three points $M,N,Q$, we conclude that $R$ is actually the center of the circle passing through $M$, $N$, and $Q$.
\end{proof}

\begin{exercise}\label{cotangentex}
Prove the Cotangent Relation in the range $\pi/4<\alpha<\pi/2$ using a similar geometric argument.
\end{exercise}

\begin{remark}
Once the Cotangent Relation is known for angles $0<\alpha<\pi/2$, then it is valid for any $\alpha$ that is not a multiple of $\pi/2$, since $\cot x$ is an odd function that is periodic with period $\pi$. If $\alpha$ were a multiple of $\pi$ then $\cot\alpha$ would be undefined; if $\alpha$ were an odd multiple of $\pi/2$, then $(\cot\alpha,0) = (0,0)$ would be on the same line as $(0,1)$ and $(0,-1)$, and hence the circle center would be undefined. This is the reason for the restriction that $\alpha$ is not a multiple of $\pi/2$.
\end{remark}

\begin{exercise}\label{right.range.ex}
Show that every \css\ contains an isosceles triangle whose vertex angle is strictly between $\pi/6$ and $5\pi/6$ in measure. (Hint: by taking a single circle center, one can find an isosceles triangle whose vertex angle is at most $2\pi/3$ in measure. What can be done if the measure of that vertex angle is too small?)
\end{exercise}

Now it is also true to say that the center of the circle passing through the three points $(0,1)$, $(0,-1)$, and $(x,0)$ is $\big((x^2-1)/2x,0\big)$. Why have we bothered to phrase the Cotangent Relation in the form we chose, instead of just stating a ``Rational Function Relation''?

The reason is that it gives us much better insight into what happens when we repeat this process of taking circle centers. Starting from the three points
\begin{equation}
M=(0,1), \quad N=(0,-1), \quad Q_1=(\cot\alpha,0),
\label{MNQ.def}
\end{equation}
let's define a sequence $\{Q_j\}$ of points recursively by setting
\begin{equation}
Q_j=C(M,N,Q_{j-1}) \quad (j\ge2).
\label{Qj.def}
\end{equation}
A simple induction using the Cotangent Relation shows that $Q_j = (\cot(2^{j-1}\alpha),0)$ for every $j\ge1$.

\begin{example}
Take $\alpha = \pi/127$, where the number $127=2^7-1$ has been carefully chosen. Then
\begin{align*}
Q_1 ={}& (\cot\frac\pi{127},0),\; Q_2 = (\cot\frac{2\pi}{127},0),\; Q_3 = (\cot\frac{4\pi}{127},0),\\
Q_4 &= (\cot\frac{8\pi}{127},0),\; Q_5 = (\cot\frac{16\pi}{127},0),\; Q_6 = (\cot\frac{32\pi}{127},0),\\
&Q_7 = (\cot\frac{64\pi}{127},0),\; Q_8 = (\cot\frac{128\pi}{127},0) = (\cot\frac\pi{127},0) = Q_1,\; Q_9 = Q_2, \dots.
\end{align*}
We see that with this starting value of $\alpha$, the sequence $\{Q_j\}$ is periodic with period~7. In fact, the same is true of the binary expansion of $\alpha/\pi = 1/127$:
\[
1/127 = (0.\ 0000001\ 0000001\ 0000001 \dots)_2.
\]
We would have needed incredibly good fortune to locate this period-7 sequence using the ``Rational Function Relation'' $(x,0) \mapsto \big( (x^2-1)/2x,0 \big)$; who would have thought to try $x=\cot(\pi/127)$?
\hfill$\diamondsuit$\end{example}

In some special situations, one of the $Q_j$ might lie at the origin, in which case $Q_{j+1} = C(M,N,Q_j)$ cannot be defined. Similarly, for some special values of $\alpha$, the evaluation $\cot(2^j\alpha)$ is undefined. It's not hard to see that these exceptions exactly correspond to each other.

\begin{example}
If we take $\alpha = \pi/32$, then we have
\begin{align*}
Q_1 ={}& (\cot\frac\pi{32},0),\; Q_2 = (\cot\frac\pi{16},0),\; Q_3 = (\cot\frac\pi8,0) = (1+\sqrt2,0),\\
Q_4 &= (\cot\frac\pi4,0) = (1,0),\; Q_5 = (\cot\frac\pi2,0) = (0,0),
\end{align*}
and so $Q_6$ cannot be defined.

The binary expansion of $\alpha/\pi = 1/32$ is simply $(0.00001)_2$, an expansion that terminates after five bits. In this case, $2^5\alpha$ is a multiple of $\pi$, and so $\cot(2^5\alpha)$ is undefined, which reflects the fact that $Q_6$ could not be defined.
\hfill$\diamondsuit$\end{example}

In general, we can see that if $Q_1 = (\cot\alpha,0)$, then $Q_j = (\cot \beta_j,0)$ where $\beta_j/\pi$ is the number whose binary expansion we get from throwing away the first $j-1$ bits in the binary expansion of $\alpha/\pi$.

\begin{exercise}
With $Q_j$ defined as in equations \eqref{MNQ.def} and \eqref{Qj.def}, show that $Q_j$ is well-defined for all $j\ge1$ unless $\alpha/\pi$ is a rational number of the form $n/2^k$. Assuming this is not the case, show that the sequence $\{Q_j\}$ is eventually periodic if and only if $\alpha/\pi$ is a rational number. For any integers $L\ge0$ and $P\ge1$, show that $\alpha$ can be chosen so that the sequence $\{Q_j\}$ is eventually periodic with period exactly $P$, after a preperiod of length exactly $L$.
\end{exercise}

\begin{example}
Consider the {\em Thu\'e-Morse constant}
\[
\tau = (0.01101001100101101001011001101001\dots)_2,
\]
where the bits in the binary expansion start with 0, then its complement 1, then the complement 10 of that pair, then the complement 1001 of all four previous bits, and so on.
%setting the first bit to be 0 and then, for each $k\ge1$, recursively setting the values of bits number $2^{k-1}+1$ to $2^k$ to be the opposites of the values of bits number 1 to $2^{k-1}$.
This builds up a sequence of truncations $\tau_k$ of $\tau$, from which the entire binary expansion of $\tau$ can be recovered:
\begin{align*}
\tau_0 &= 0.\ 0, \\
\tau_1 &= 0.\ 0\ 1, \\
\tau_2 &= 0.\ 01\ 10, \\
\tau_3 &= 0.\ 0110\ 1001, \\
\tau_4 &= 0.\ 01101001\ 10010110,\, \dots.
\end{align*}
It is not hard to prove that there are never three 0's or three 1's in a row in this binary expansion, which has an interesting implication (see Exercise \ref{TM.exercise}) for the sequence of points $Q_j$ we get by taking $\alpha=\tau/\pi$ in equations \eqref{MNQ.def} and \eqref{Qj.def}. One important to fact to note is that the binary expansion of $\tau$ is not eventually periodic, and hence all of the $Q_j$ are distinct in this case. (See \cite{EW} for other properties of this interesting number.)
\hfill$\diamondsuit$\end{example}

\begin{exercise}
Prove that the Thu\'e-Morse constant $\tau$ never has three 0's in a row or three 1's in a row in its binary expansion. If $Q_j$ is defined as in equations \eqref{MNQ.def} and \eqref{Qj.def} with the choice $\alpha=\tau/\pi$, show that all the $Q_j$ are on the line segment between $(\cot(7\pi/8),0) = (-1-\sqrt2,0)$ and $(\cot(\pi/8),0) = (1+\sqrt2,0)$.
\label{TM.exercise}
\end{exercise}

From these examples, it seems that lots of choices for $\alpha$ result in the sequence $\{Q_j\}$ being bounded. As it turns out, this is quite misleading: while such examples are easy to write down, they are actually relatively rare. It is known that a ``typical'' real number has every possible string of bits in its binary expansion, including a string of 0's or 1's as long as we wish. In fact, even more is true: if we count how many strings of length $k$ in the binary expansion of a ``typical'' real number are equal to our favorite length-$k$ string, then the result becomes closer and closer to $1/2^k$ as we go further and further out. This property of a real number is called being {\em normal to the base~2}. For example, the pairs of bits in the binary expansion of a ``typical'' real number are, in the limit, evenly split among the $4=2^2$ possibilities 00, 01, 10, and 11.

So what do we mean by ``typical''? Suppose a set $S$ has the property that for every positive number $\ep$, there is a collection of open intervals (maybe infinitely many intervals) of total length at most $\ep$ whose union contains $S$. Then in many ways $S$ acts like a ``very small'' set. Any countable set, such as $\Q$, is a set of this type: simply choose an open interval of length $\ep/2$ containing the first element of the set, an open interval of length $\ep/4$ containing the second element, an open interval of length $\ep/8$ containing the third element, and so on. The technical term is that $S$ ``has measure zero'', in the terminology of Lebesgue measure (see \cite[Chapter 3]{HR}), and the idea is that $S$ is such a small set that numbers in $S$ are ``atypical'', while numbers outside of $S$ are ``typical''. In fact we say that ``almost all'' numbers are outside of $S$.

One of the fundamental theorems in this subject is that the set of real numbers that are not normal to the base~2 has measure zero. This doesn't mean that every real number is normal to the base~2---in fact, none of the rational numbers are. (See the article \cite{GM} for a construction of a very abnormal number in this respect.) What this result does mean is that ``almost all'' real numbers are normal to the base~2: if we chose one at random, we would have a ``100\% chance'' of picking a normal one!

Since almost all real numbers are normal, that means that almost all real numbers contain any given binary string somewhere (indeed, many times) in their binary expansion. Therefore by throwing away the right number of initial bits, that given binary string can be brought to the front of the binary expansion, meaning that the resulting number is as close to any prescribed real number as we like. The following corollary of these observations gives the big picture for this one-dimensional problem:

\begin{theorem}
Let $M=(0,1)$ and $N=(0,-1)$. For almost all choices of a point $Q_1$ on the $x$-axis, the closure of the sequence $\{Q_j\}$ defined by equations \eqref{MNQ.def} and \eqref{Qj.def} is the entire $x$-axis.
\end{theorem}

\section{The plane truth}

Now let's go back to two dimensions and investigate Theorems \ref{dense.thm} and ~\ref{only1.thm}. We start by proving the interesting fact that two triangles that are naturally related in the context of taking circle centers are actually similar to each other. This Similarity Relation is a more developed version of the Cotangent Relation we saw earlier (and so our detour into one dimension really was helpful!).

\newtheorem*{SRsc}{Similarity Relation (special case)}
\begin{SRsc}
Let $\alpha$ be a real number in the range $0<\alpha<\pi/2$. Let $M=(0,1)$, $N=(0,-1)$, $Q=(\cot\alpha,0)$, and $R=(\cot2\alpha,0)$, and define $S$ to be the center of the circle passing through $M$, $Q$, and $R$. Then $\triangle NMQ$ is similar to $\triangle QRS$, both triangles being isosceles with vertices $Q$ and $S$, respectively. The constant of proportionality is $(\cot\alpha-\cot2\alpha)/2$, and the bases of $\triangle NMQ$ and $\triangle QRS$ are perpendicular.
\end{SRsc}

We note that the function $\cot x$ is strictly decreasing on the interval $0<x<\pi$, and so $(\cot\alpha-\cot2\alpha)/2$ is always positive in the range $0<\alpha<\pi/2$.

\begin{proof}
We know from the Cotangent Relation that $R$ is actually $C(M,N,Q)$ in disguise. This implies that $MR=QR$, and so $\triangle MQR$ is isosceles with vertex $R$. Now it is an elementary fact from Euclidean geometry that if $S$ is defined to be the center of the circle through $M$, $Q$, and $R$, then $\overline{RS}$ is perpendicular to $\overline{MQ}$. Indeed, both $R$ and $S$ lie on the perpendicular bisector of $\overline{MQ}$, the former because $\triangle MQR$ is isosceles. As we already know that the measure of $\angle MQR$ is $\alpha$, this implies that the measure of $\angle QRS$ is $\pi/2-\alpha$, as Figure~\ref{SRfig}(a) shows. We now see that $\triangle MNQ$ and $\triangle QRS$ are both isosceles triangles with base angles of $\pi/2-\alpha$ (Figure~\ref{SRfig}(b)), and therefore these two triangles are similar.

%:SRfig
\begin{figure}[hbtp]

\hfill
\begin{picture}(0, 115)
\put(-25,100){(a)}
\put(-10, 48){$M$}
\put(91, 25){$Q$}
\put(36, 16){$R$}
\put(63, 110){$S$}
\put(64, 29){$\scriptscriptstyle\alpha$}
\put(72, 55){\vector(-1,-1){25}}
\put(73, 57){$\scriptstyle\pi/2-\alpha$}
\end{picture}
\includegraphics[height=1.5in]{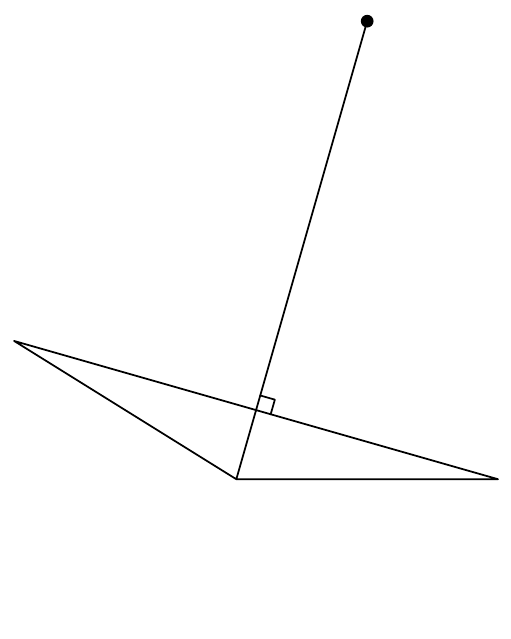}
\hfill\hfill
\begin{picture}(0,0)
\put(-25,100){(b)}
\put(-10, 48){$M$}
\put(-8, -1){$N$}
\put(91, 25){$Q$}
\put(33, 25){$R$}
\put(63, 109){$S$}
\put(16, 66){$\scriptstyle\pi/2-\alpha$}
\put(32, 62){\vector(1,-2){16}}
\put(28, 62){\vector(-4,-3){20}}
\end{picture}
\includegraphics[height=1.5in]{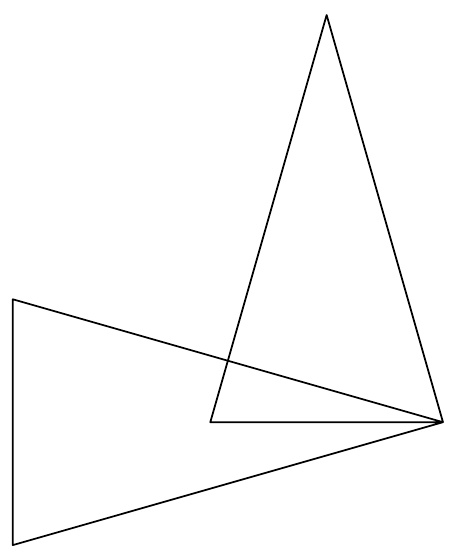}
\hfill\hfill

\caption{The Similarity Relation}
\label{SRfig}
\end{figure}
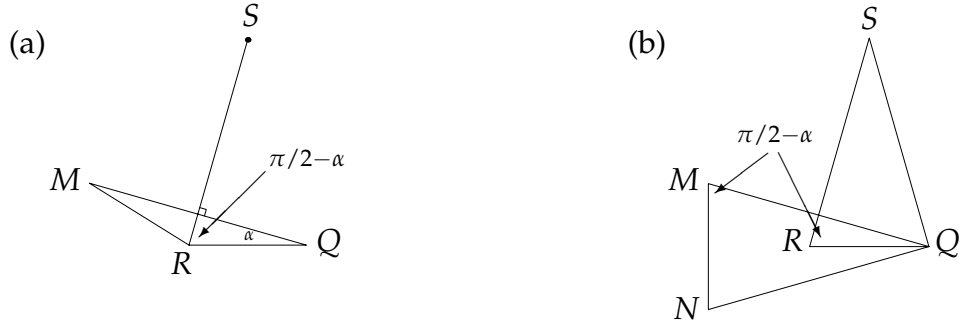

The fact that the constant of proportionality equals $(\cot\alpha-\cot2\alpha)/2$ follows from noting that $MN=2$ while $QR=\cot\alpha-\cot2\alpha$. Just as easily, we see that the base $\overline{QR}$ of $\triangle QRS$ lies upon the perpendicular bisector of the base $\overline{MN}$ of $\triangle MNQ$, by the same argument as in the previous paragraph. In particular, the two triangles' bases $\overline{MN}$ and $\overline{QR}$ are indeed perpendicular.
\end{proof}

Again, we can scale and rotate the plane without affecting these relative statements about circle centers. Therefore, after the slight change of notation $\beta=2\alpha$, we have the following corollary:

\newtheorem*{SR}{Similarity Relation}
\begin{SR}
Let $\triangle DEF$ be isosceles with vertex $F$, and let $\beta$ be the measure of the vertex angle at $F$. Let $G=C(D,E,F)$ and $H=C(E,F,G)$. Then $\triangle FGH$ is similar to $\triangle DEF$ with constant of proportionality
\begin{equation}
\lambda = (\cot(\beta/2)-\cot\beta)/2.  \label{lambdadef}
\end{equation}
Moreover, the axes of symmetry of the two triangles are perpendicular.
\end{SR}

We remark that the constant $\lambda$ defined in equation \eqref{lambdadef} is less than 1 precisely when $\beta$ is in the range $\pi/6 < \beta < 5\pi/6$.

Now let's play a game similar to the one we played last section---namely, recursively generate new points by taking centers of circles going through old points. This time, however, instead of holding two of the points fixed, we let all of the new points take their fair turns.

As remarked earlier, every \css\ contains an isosceles triangle, so let's start with an isosceles $\triangle P_1P_2P_3$ with vertex $P_3$, and let $\beta$ denote the measure of the vertex angle:
\begin{equation}
P_3P_1 = P_3P_2, \quad \beta = \mathop{\text{m}}(\angle P_1P_3P_2).
\label{P123.def}
\end{equation}
From these starting points, we recursively define a sequence of circle-centers:
\begin{equation}
P_n=C(P_{n-1},P_{n-2},P_{n-3}) \quad (n\ge4).
\label{Pn.def}
\end{equation}
By the Similarity Relation, we know that $\triangle P_3P_4P_5$ is similar to $\triangle P_1P_2P_3$, where the constant of proportionality $\lambda$ is given in equation \eqref{lambdadef}. Furthermore, $\triangle P_3P_4P_5$ is itself isosceles with vertex angle $\beta$, and so $\triangle P_5P_6P_7$ is similar to $\triangle P_3P_4P_5$ by the Similarity Relation again, with the same constant of proportionality. Continuing in this way as long as we wish, we see that for any number $k\ge1$, the triangle $\triangle P_{2k+1}P_{2k+2}P_{2k+3}$ is similar to $\triangle P_1P_2P_3$ with constant of proportionality $\lambda^k$.

\begin{figure}[hbtp]

\hfill
\begin{picture}(0, 0)
\put(-25,100){(a)}
\put(97, 71){$\scriptscriptstyle\infty$}
\end{picture}
\includegraphics[height=2.25in]{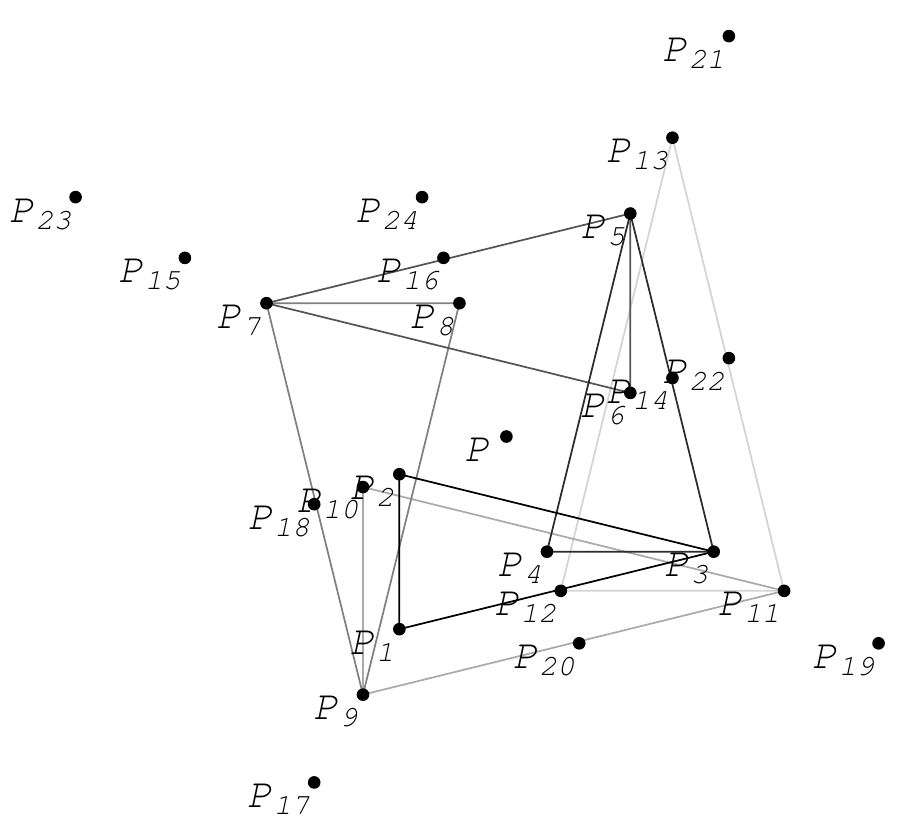}
\hfill\hfill
\begin{picture}(0,0)
\put(-25,100){(b)}
\put(72, 82){$\scriptscriptstyle\infty$}
\end{picture}
\includegraphics[height=2.25in]{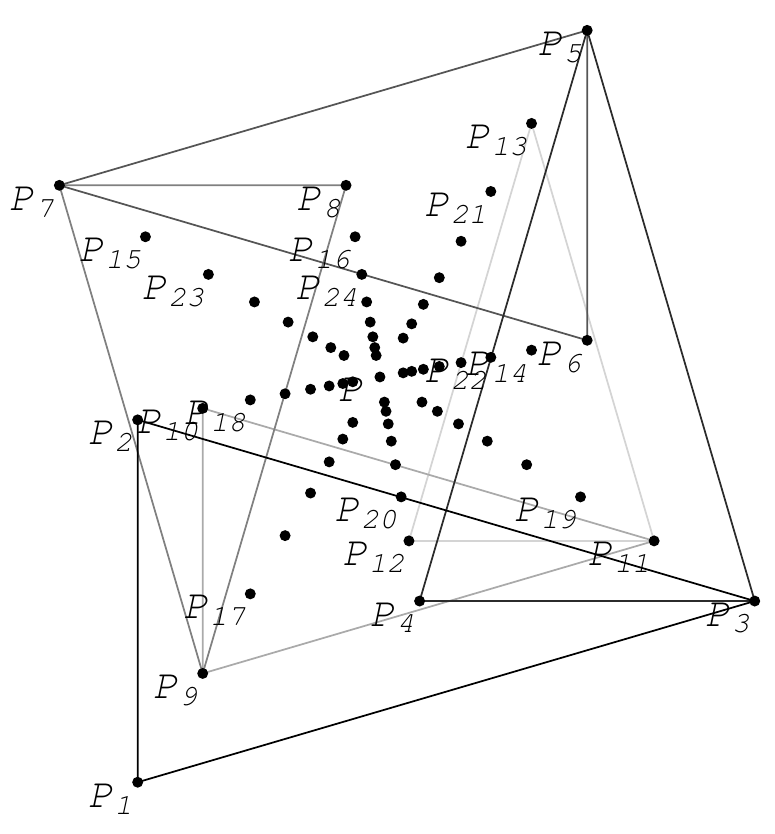}
\hfill\hfill

\caption{Points generated by recursively taking circle centers}
\label{recursefig}
\end{figure}
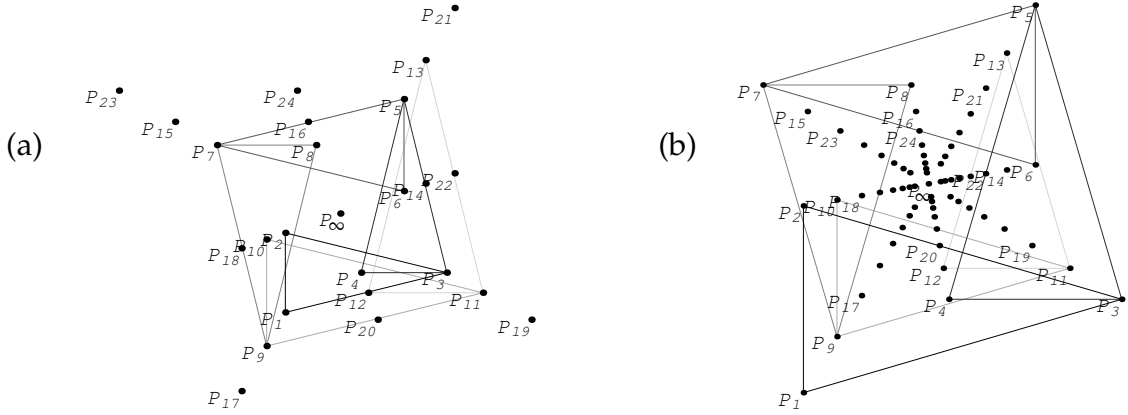

The result is that these triangles swirl around in the plane, getting bigger or smaller depending on whether or not $\lambda>1$. Figure~\ref{recursefig}(a) was drawn with $\beta=2\pi/13$, corresponding to $\lambda\approx1.076$, while Figure~\ref{recursefig}(b) was drawn with $\beta=2\pi/11$, corresponding to $\lambda\approx0.925$. In both figures, the darkest triangle is the original $\triangle P_1P_2P_3$, the successive triangles getting lighter and lighter the farther along we go in the recursion.

%:recursefig

Spurred on by the prominent patterns in Figure~\ref{recursefig}, let's examine in more detail the behavior of the points $P_n$. For notational purposes, we define the following sets of points:
\begin{gather*}
\S_1 = \{P_1,P_5,P_9,P_{13},\dots\}, \quad \S_2 = \{P_2,P_6,P_{10},P_{14},\dots\}, \\
\S_3 = \{P_3,P_7,P_{11},P_{15},\dots\}, \quad \S_4 = \{P_4,P_8,P_{12},P_{16},\dots\},
\end{gather*}
where in each set the indices increase by four.

\newtheorem*{OQ}{Orderly Queues Relation}
\begin{OQ}
Let $P_1,P_2,P_3$ be points in the plane such that $\triangle P_1P_2P_3$ is isosceles with vertex $P_3$. Define the sequence $P_n$ recursively by $P_n=C(P_{n-1},P_{n-2},P_{n-3})$ for $n\ge4$. Then there exists a point $P_\infty$ in the plane with the following properties:
\begin{enumerate}
\item Each of the sets $\S_1$, $\S_2$, $\S_3$, and $\S_4$ is contained in a line that goes through $P_\infty$;
\item The lines containing $\S_1$ and $\S_3$ are perpendicuar to each other, as are the lines containing $\S_2$ and $\S_4$.
\end{enumerate}
\end{OQ}

These orderly queues can actually be seen in Figure~\ref{recursefig}, even in Figure~\ref{recursefig}(a) where the points $P_n$ are heading away from $P_\infty$, but certainly in Figure~\ref{recursefig}(b) where they are heading towards~$P_\infty$.

\begin{proof}
Out of all compositions of translations, scalings, and rotations of the plane, there is a unique mapping $f$ that sends $P_1$ to $P_3$ and $P_2$ to $P_4$; this mapping has a unique fixed point (see Exercise \ref{mapping.exercise} below), which we label $P_\infty$, and the mapping $f$ can be thought of as the composition of a scaling and a rotation, both centered at $P_\infty$. As $\vvec{P_2P_1}$ and $\vvec{P_4P_3}$ are perpendicular, the rotation component must be a rotation by $\pi/2$. If we were to set coordinates so that $P_\infty$ were the origin, then $f$ could be represented simply as multiplication by a $2\times2$ matrix:
\[
f\bigg( \bigg( \begin{matrix} x \\ y \end{matrix} \bigg) \bigg) = \lambda \bigg( \begin{matrix} \cos\pi/2 & {-\sin\pi/2} \\ \sin\pi/2 & \cos\pi/2 \end{matrix} \bigg) \bigg( \begin{matrix} x \\ y \end{matrix} \bigg) = \bigg( \begin{matrix} 0 & \hskip-6pt {-\lambda} \\ \lambda & \hskip-6pt0 \end{matrix} \bigg) \bigg( \begin{matrix} x \\ y \end{matrix} \bigg),
\]
with $\lambda$ as defined in equation~\eqref{lambdadef}.

Now the Similarity Relation can be interpreted as saying that the mapping $f$ also sends $P_3$ to $P_5$. Moreover, as we have mentioned before, translations, scalings, and rotations all respect circumcenters of triangles (that is, under any of these transformations, the image of a triangle's circumcenter is the same as the circumcenter of the triangle's image). Therefore we have
\[
f(P_4) = f(C(P_1,P_2,P_3)) = C(f(P_1), f(P_2), f(P_3)) = C(P_3,P_4,P_5) = P_6,
\]
and an easy induction using the same argument shows that in fact $f(P_n) = P_{n+2}$ for every $n\ge1$. Since $f$ causes a rotation by $\pi/2$ around the point $P_\infty$, the second iterate $f\circ f$ causes a rotation by $\pi$, and so the points $P$, $P_\infty$, and $f(f(P))$ always lie a line. However, $f(f(P_n)) = P_{n+4}$ for every $n\ge1$. Therefore all the points in any of the sets $\S_1,\S_2,\S_3,\S_4$ lie on a single line that includes $P_\infty$, as claimed. The argument also shows that the lines containing $\S_1$ and $\S_3$ are perpendicuar to each other, as are the lines containing $\S_2$ and~$\S_4$. 
\end{proof}

\begin{remark}
If we choose the first three points to be $(0,-1)$, $(0,1)$, and $(0,x)$, for instance, then it can be worked out by the compulsive reader that the coordinates of $P_\infty$ are
\[
P_\infty = \bigg( \frac{12x^3-4x}{x^4+18x^2+1} , \frac{3x^4+2x^2-1}{x^4+18x^2+1} \bigg).
\]
\end{remark}

\begin{exercise}
Let $P_1$, $P_2$, $Q_1$, and $Q_2$ be points in the plane, with $P_1\ne P_2$ and $Q_1\ne Q_2$. Show that among all compositions of translations, scalings, and rotations, there is a unique mapping that sends $P_1$ to $Q_1$ and $P_2$ to $Q_2$. Moreover, if the vectors $\vvec{P_1P_2}$ and $\vvec{Q_1Q_2}$ are distinct, show that this mapping has a unique fixed point.
\label{mapping.exercise}
\end{exercise}

\begin{remark}
This exercise can be done using conventional Euclidean geometry; however, it is instructive to find a solution by viewing the plane as the set of complex numbers, whereupon the mappings described are simply the linear functions $z\mapsto az+b$.
\end{remark}

The Orderly Queues Relation provides us with quite a bit of structure hidden inside any \css. Is there a way to exploit this structure to prove that every \css\ must be dense in the plane---that is, that any closed \css\ must fill up the entire plane down to the last point? The key step on the path to answering this question turns out to be showing that any closed \css\ must at least fill up an entire line segment.

\newtheorem*{SFF}{Segment-Filling Fact}
\begin{SFF}
Let $\S$ be a closed \css, and suppose that there are two perpendicular lines $\L_1$ and $\L_2$ such that $\S$ contains:
\begin{itemize}
\item the point of intersection $I$ of $\L_1$ and $\L_2$;
\item another point $A$ on $\L_1$ besides $I$; and
\item a sequence of points \{$P_1$, $P_2$, \dots\} on $\L_2$ that converges to $I$.
\end{itemize}
Then $\S$ contains the entire line segment $\overline{AI}$.
\end{SFF}

\begin{proof}
As usual, we may scale, rotate, and translate our set $\S$ as we wish, and so without loss of generality we can assume the following: $\L_1$ is the $y$-axis, $\L_2$ is the $x$-axis, $I$ is the origin, and $A=(0,1)$. The sequence of points \{$P_1$, $P_2$, \dots\} then lies along the $x$-axis. To prove that $\S$ contains the entire segment $\overline{AI}$, it suffices to prove that $\S$ contains every point of the form $(0,a/2^n)$ where $n\ge0$ and $0\le a\le 2^n$ are integers, since the closure of this set of points is $\overline{AI}$.

In fact, we choose to prove an even stronger statement: for every integer $n\ge0$, and for every integer $0\le a\le 2^n$, not only does $\S$ contain the point $(0,a/2^n)$, but it also (except in the case $a/2^n = 1$) contains a sequence of points on the line $y=a/2^n$ that converges to $(0,a/2^n)$. (Figure \ref{inductionfig} helps us visualize just what is being claimed here. In that figure, the $x$-coordinates of the $\{P_j\}$ are shown decreasing monotonically to 0; this could be assumed anyway without losing generality by passing to a subsequence, but it isn't important to the proof.) The reason we make our job apparently harder in this way is that the new statement can actually be proved by induction on $n$.

%:inductionfig
\begin{figure}[hbtp]

\hfill
\begin{picture}(0, 0)
\put(-25,150){(a)}
\end{picture}
\includegraphics[height=2.6in]{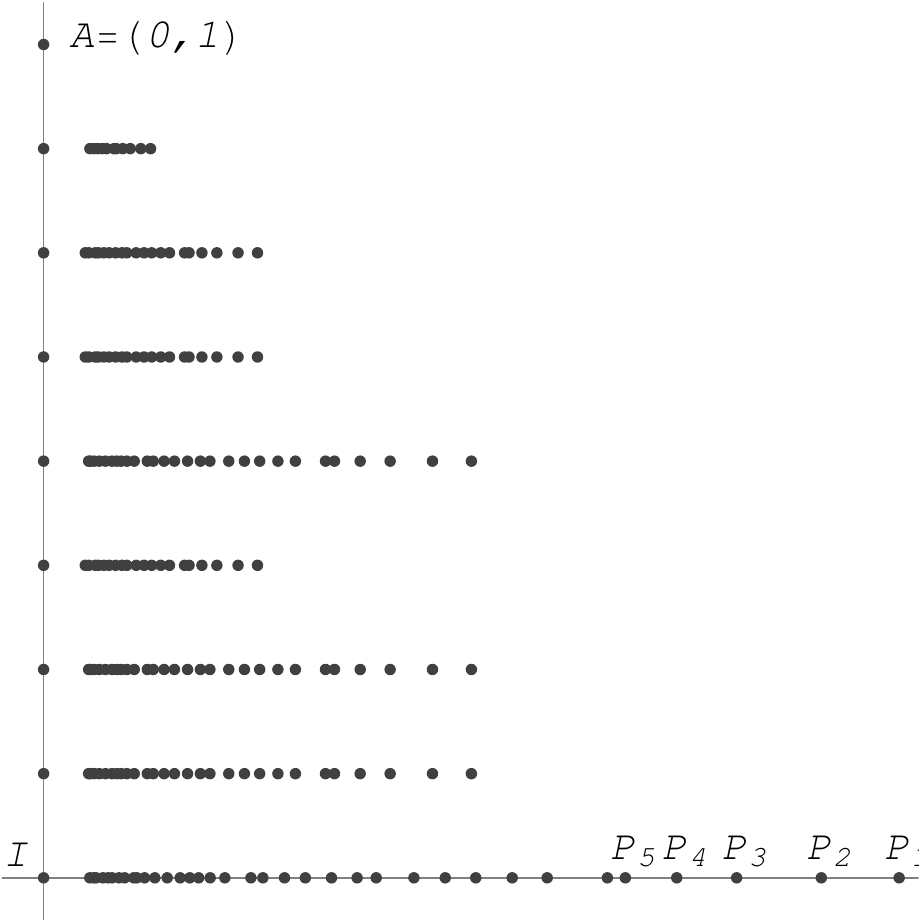}
\hfill\hfill
\begin{picture}(0,0)
\put(0,150){(b)}
\end{picture}
\includegraphics[height=2.6in]{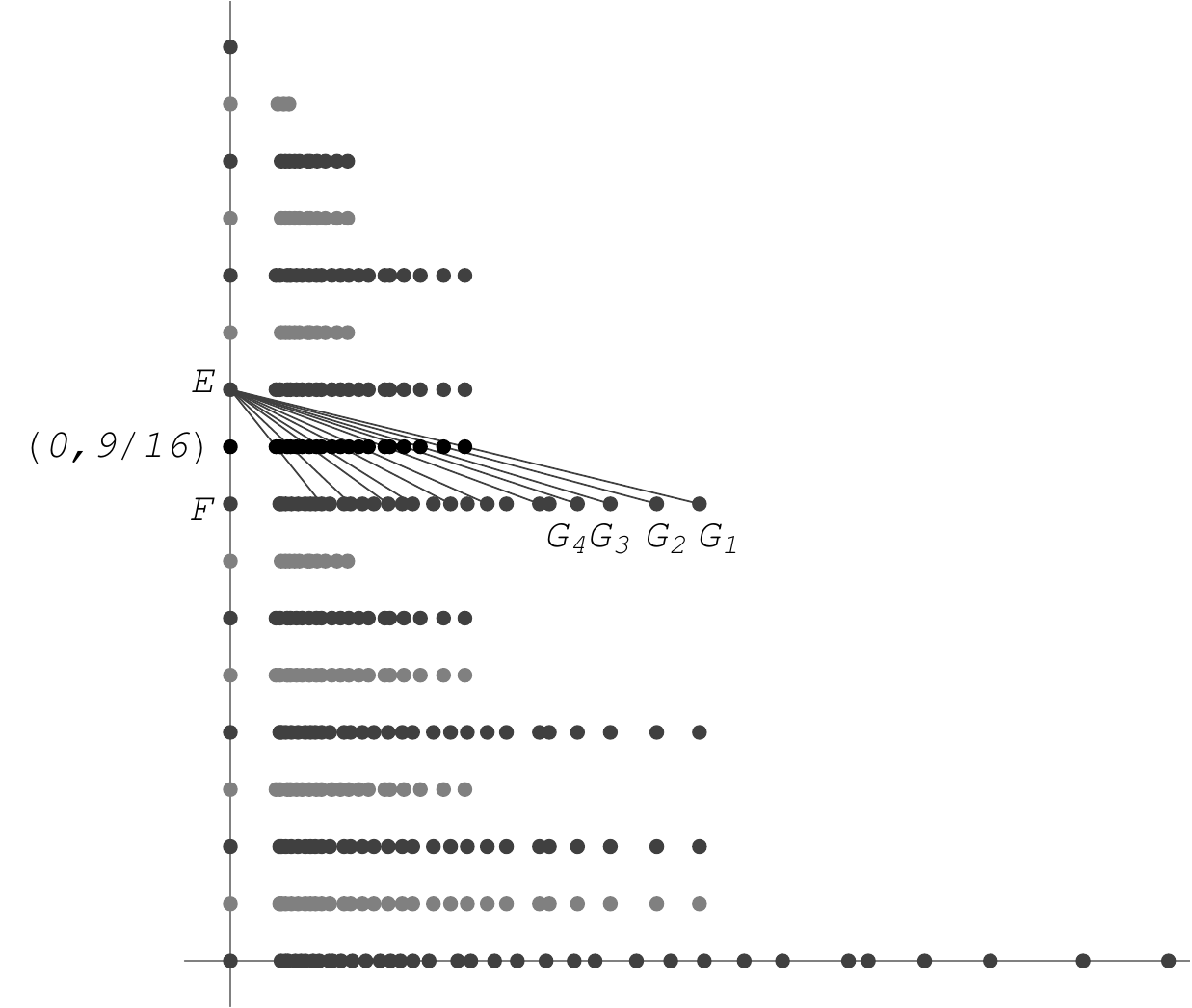}
\hfill\hfill

\caption{Inductive proof of the Segment-Filling Fact}
\label{inductionfig}
\end{figure}
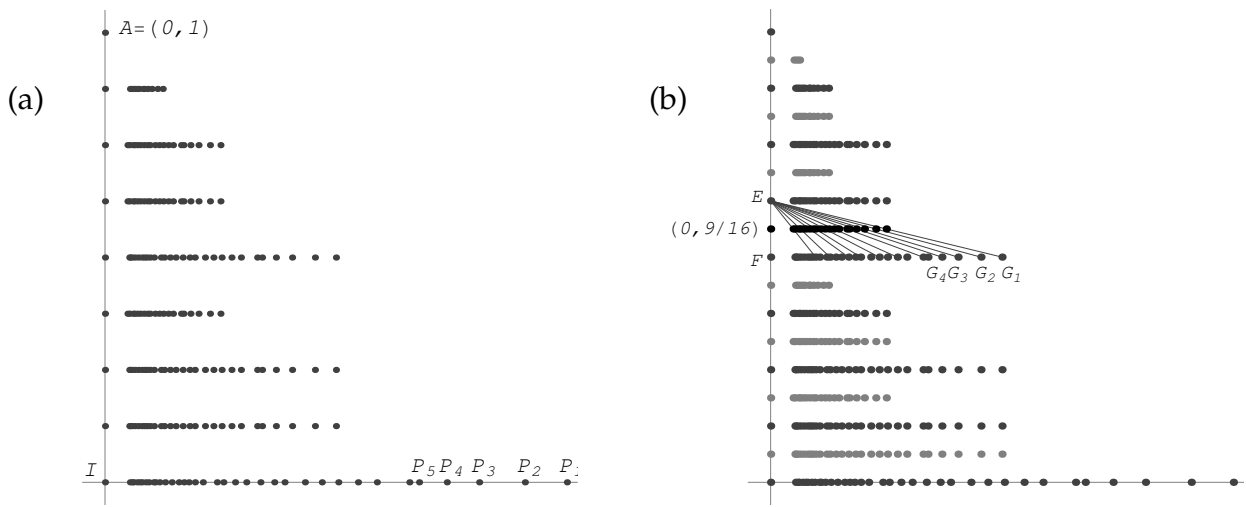

First of all, note that the base case $n=0$ is true precisely because of the hypotheses described in the first paragraph of this proof. For the inductive step, assume that the case $n-1$ is known to be true (Figure \ref{inductionfig}(a) provides a picture of this induction hypothesis for $n-1=3$); we need to establish that the case $n$ is true as well. If the numerator $a$ is even, say $a=2b$, then $a/2^n=b/2^{n-1}$, and so the induction hypothesis takes care of those values of $a$ immediately. Thus we can restrict our attention to odd numerators $a=2b-1$ with $1\le b\le 2^{n-1}$. Note that it suffices to prove that $\S$ contains a sequence of points on the line $y=(2b-1)/2^n$ that converges to $(0,(2b-1)/2^n)$: since $\S$ is closed, it will then automatically contain the point $(0,(2b-1)/2^n)$ as well.

Now by the induction hypothesis, $\S$ contains both $E=(0,b/2^{n-1})$ and $F=(0,(b-1)/2^{n-1})$, as well as a sequence of points $\{G_j\}$ on the line $y=(b-1)/2^{n-1}$ converging to $F$. (Figure \ref{inductionfig}(b) illustrates this with $n=4$ again and $b=5$.) Note that $\triangle EFG_j$ is a right triangle with hypotenuse $\overline{EG_j}$. Since the circumcenter of a right triangle is the midpont of its hypotenuse, the \css\ $\S$ also contains the midpoint of each of the segments $\overline{EG_j}$. But it is easily seen that this sequence of midpoints lies on the line $y=\tfrac12((b-1)/2^n + b/2^n) = a/2^n$ and converges to $(0,a/2^n)$. This establishes the inductive step and hence finishes the proof.
\end{proof}

\begin{exercise}
If a closed \css\ $\S$ contains a sequence of points on a line $\L$ that converges to a point $P$, prove that it automatically contains a point on the line through $P$ that is perpendicular to $\L$.
\end{exercise}

Since we're obtaining more and more information about what's happening inside \css s, let's press on to another useful fact:

\newtheorem*{QF}{Quadrilateral Fact}
\begin{QF}
Any \css that contains an entire line segment must also contain an entire quadrilateral and its interior.
\end{QF}

(By ``quadrilateral'' we mean a nondegenerate quadrilateral, of course, so that it has an interior to speak of.)

\begin{proof}
By the usual rotation/translation/scaling argument, we can assume that our \css\ contains a segment of the $x$-axis together with the point $(0,1)$. It then remains to understand the behavior of the function
\[
f(u,v) = \text{the center of the circle through $(0,1)$, $(u,0)$, and $(v,0)$} = \bigg( \frac{u+v}{2},\frac{u v+1}{2} \bigg)
\]
which takes values in the plane. The following exercise does the rest of the work.
\end{proof}

\begin{exercise}\label{quadex}
Let $a<b<c<d$ be real numbers. Show that as $u$ ranges over the interval $[a,b]$ and $v$ ranges over the interval $[c,d]$, the function $f(u,v)$ fills the quadrilateral whose vertices are $f(a,c)$, $f(a,d)$, $f(b,d)$, and $f(b,c)$. Conclude that the Quadrilateral Fact is true in general.
\end{exercise}

Finally, we sprint to the end of our journey by assembling these known facts together into a proof of our two (equivalent) main theorems.

\begin{proof}[Proof of Theorem \ref{only1.thm}]
Let $\S$ be a closed \css; we want to show that $\S=\R^2$. By Exercise \ref{right.range.ex}, we know that $\S$ contains an isosceles triangle with vertex angle strictly between $\pi/6$ and $5\pi/6$ in measure. Because the constant $\lambda$ defined in equation~\ref{lambdadef}  is less than 1 for this vertex angle, the Orderly Queues Relation tells us that $\S$ contains a sequence of points $\S_3$ which all lie on a common line and which converge to the point $P_\infty$. In this case, the point $P_\infty$ is also in $\S$ since $\S$ is closed. Moreover, $\S$ contains the point $P_1$, which is on the line through $P_\infty$ that is perpendicular to the line containing $\S_3$. Therefore, by the Segment-Filling Fact, $\S$ must contain the entire line segment joining $P_\infty$ and $P_1$.

Once the \css\ $\S$ contains a line segment, the Quadrilateral Fact shows that $\S$ must actually contain an entire quadrilateral and its interior. From here we can quickly show that $\S$ must be all of $\R^2$, by the following argument: Choose any point $A$ in $\R^2$, and consider the circle centered at $A$ that goes through the center of the quadrilateral (any circle centered at $A$ that intersects the interior of the quadrilateral would do). The intersection of this circle with the interior of the quadrilateral is some circular arc. No matter how small this arc might be, we can choose three distinct points on the arc, which are all in $\S$ because they're inside the quadrilateral. Therefore the center of the circle through those three points must be in $\S$---but that center is just the point $A$ we chose! And since $A$ was an arbitrary point in the plane, we see that $\S$ must be all of $\R^2$, as claimed.
\end{proof}

{\it Acknowledgements.} The author thanks Lowell W.\ Beineke and an anonymous referee for helpful comments on an earlier version of this manuscript; several of their suggestions have become improvements in the exposition of this paper. The author was supported in part by grants from the Natural Sciences and Engineering Research Council.

\end{document}